\documentclass[12pt,reqno,oneside]{amsart}

\usepackage{xcolor}
\usepackage{amsmath,amsthm}
\usepackage{graphicx}
\usepackage{amssymb}%??????
\usepackage{mathrsfs}%????
\newcommand\Wtilde{\stackrel{\sim}{\smash{\mathcal{M}}\rule{0pt}{1.2ex}}}

\newtheorem{thm}{Theorem}[section]

\newtheorem{defn}[thm]{Definition}

\newtheorem{rmk}[thm]{Remark}
\newtheorem{nt}[thm]{Notation}
\newcommand*\diff{\mathop{}\!\mathrm{d}}
\newcommand*{\ccdot}{\raisebox{-0.25ex}{\scalebox{1.2}{$\cdot$}}}
\makeatletter

\newcommand{\Rmnum}[1]{\expandafter\@slowromancap\romannumeral #1@}
\makeatother
\newcommand{\overbar}[1]{\mkern 1.5mu\overline{\mkern-1.5mu#1\mkern-1.5mu}\mkern 1.5mu}

\newcommand\blfootnote[1]{%
  \begingroup
  \renewcommand\thefootnote{}\footnote{#1}%
  \addtocounter{footnote}{-1}%
  \endgroup
}
\newcommand{\Addresses}{{% additional braces for segregating \footnotesize
  \bigskip
  \footnotesize

\textsc{College of Liberal Arts and Sciences, National University of Defense Technology,  Kaifu District, 
Changsha, China}\par\nopagebreak
  \textit{E-mail address}: \texttt{zhouzijian.edu@gmail.com}

}}
\begin{document}
\title[Ekedahl-Oort strata of curves of genus four in char. three]{Ekedahl-Oort strata of curves of genus four in characteristic three}
\author{Zijian Zhou}
\maketitle
\vspace{-1cm}
\begin{abstract}
We study the induced Ekedahl-Oort stratification on the moduli space $\mathcal{M}_4$ of curves  of genus $4$ in characteristic $3$. By constructing families of curves with given Ekedahl-Oort type and by computing the dimension of the intersection of the induced Ekedahl-Oort strata and the boundary divisor classes, we show that for certain induced Ekedahl-Oort strata in $\mathcal{M}_4$, they have the same codimension in $\mathcal{M}_4$ as the corresponding Ekedahl-Oort strata in moduli space $\mathcal{A}_4$ of principally polarized abelian varieties of dimension $4$.
\end{abstract}
\thispagestyle{empty}
\vspace{-0.6cm}
\section{Introduction}
Let $k$ be an algebraically closed field of characteristic $p>0$. 
Let $\mathcal{A}_g$ be the moduli space (stack) of principally polarized
abelian varieties of dimension $g$ defined over $k$ and let $\mathcal{M}_g$ be the
moduli space of (smooth projective) curves of genus~$g$ defined over $k$.
Ekedahl and Oort introduced a stratification on $\mathcal{A}_g$ consisting of $2^g$ strata, cf. \cite{Ekedahl2009,Oort1999}. It characterizes the $p$-torsion group scheme of abelian varieties. Later Ekedahl and van der Geer \cite{Ekedahl2009,vanderGeer1999} realized this stratification can be defined using the de Rham cohomology.  These strata are indexed by $n$-tuples
$\mu=[\mu_1, \ldots,\mu_n]$ with $0\leq n \leq g$ and $\mu_1 > \mu_2
>\cdots >\mu_n>0$. The largest stratum is the locus
of ordinary abelian varieties corresponding to the empty $n$-tuple $\mu=\emptyset$. There are many results about the Ekedahl-Oort stratification on $\mathcal{A}_g$. For example, Oort \cite{Oort1999} showed that  any positive dimensional stratum is connected and the Ekedahl-Oort stratum index by $\mu$ is of codimension $\sum_{i=1}^n\mu_i$ in $\mathcal{A}_g$. Ekedahl and van der Geer \cite{vanderGeer1999} studied the irreducibility of certain Ekedahl-Oort strata and computed  the Chern classes defined by Ekedahl-Oort strata .

Via the Torelli map $\tau:
\mathcal{M}_g\to \mathcal{A}_g$  this stratification can be pulled back
to $\mathcal{M}_g$ and one can  ask what this
stratification is.  The curves of genus $4$ in characteristic $3$ are here of special interest.  For hyperelliptic curves with $p = 2$, Elkin and Pries  \cite{MR3095219} describe  their Ekedahl-Oort types completely. In \cite{2018arXiv181204996Z}, we determine the dimension and reducibility of Ekedahl-Oort strata of hyperelliptic curves in characteristic $3$.

 In this paper, we study the dimension of  Ekedahl-Oort strata in $\mathcal{M}_4$. For certain Ekedahl-Oort strata, we determine the dimension by looking at the intersection of the Ekedahl-Oort strata and the boundary divisor class. We also show the codimension of Ekedahl-Oort strata satifying certain bounds by constructing families of curves lies in given Ekedahl-Oort strata.    One can find related results in \cite{2017arXiv171204921D}. 
 
 To describe the results precisely, we need the following notation.
 \begin{nt}
 For the Ekedahl-Oort stratum indexed by $\mu$, we denote the induced strata on $\mathcal{M}_g$ by $Z_{\mu}$. Note that  the indices $\mu$ of the Ekedahl-Oort strata are partially ordered by  $$\mu=[\mu_1,\dots,\mu_n] \preceq \upsilon=[\upsilon_1,\dots,\upsilon_m]$$
if $n\leq m$ and $\mu_i \leq\upsilon_i $ for $i=1,\dots,n$.
 We say a (smooth) curve has Ekedahl-Oort type $\mu$ if the corresponding point in $\mathcal{M}_g$ lies in $Z_{\mu}$.  
 \end{nt}

\begin{thm}\label{thmmain}
Let $p=3$. The locus $Z_{\mu}$ is non-empty of codimension $\sum_{i=1}^n\mu_i$ in $\mathcal{M}_4$ if $ \mu  \preceq [4,1]$. Moreover,  the locus $Z_{\mu}$ is non-empty of codimension at most $ \sum_{i=1}^n\mu_i$ for $\mu=[3,2], [3,2,1]$ and $Z_{[4,3]}$ is non-empty of dimension $3$.
\end{thm}
Part of Theorem \ref{thmmain}  was known. Pries \cite[Theorem 4.3]{MR2569747} showed that there is a family (with dimension $6$) of smooth curves of genus $4$ with $p$-rank 1 and $a$-number 1(that is , curve with Ekedahl-Oort type $\mu=[3]$). Moreover, if $p\geq 5$, then by \cite[Corollary 4.5]{MR2569747}   there is a family (with dimension $6$) of smooth curves of genus $4$ with $p$-rank 2 and $a$-number 2 (Ekedahl-Oort type $\mu=[2,1]$).

\section{Background}
From now on , let $k$ be an algebraically closed field of characteristic $3$. All objects are defined on the category of $k$-schemes.
\subsection{The first de Rham cohomology}
Suppose Y is a scheme. Here we recall the first de Rham cohomology of $Y$ using the C\v ech cohomology. Denote by $(\Omega_Y^{\ccdot},\diff )$ the complex of sheaves of K\"ahler differential forms on $Y$. The  de Rham cohomology of $Y$ is  defined as the hypercohomology of the functor $H^0(Y,\ccdot)$ with respect to $(\Omega_Y^{\ccdot},\diff )$, cf. \cite[Section 11.4]{EGA}. Since we mainly focus on curves and abelian varieties, we are only interested in the first de Rham cohomology, which by  \cite{MR0241435} can also be described in terms of C\v ech cocycles as follows.
Let $\mathcal{U}=\{U_i | i\in I\}$ be an  open affine cover of $Y$. For  a sheaf  $\mathscr{F}$ of abelian groups,  we denote the C\v ech complex of abelian groups associated to $\mathscr{F}$   by  $C^{\ccdot}(\mathcal{U},\mathscr{F})$  \cite[Chapter  \Rmnum{3}, Section 4]{MR0463157}. We define the $H^1_{dR}(\mathcal{U})$ with respect to the covering $\mathcal{U}$ to be
\begin{align}
H^1_{dR}(\mathcal{U})=Z_{dR}^1(\mathcal{U})/B_{dR}^1(\mathcal{U}), \label{de Rham defi}
\end{align}
where $Z_{dR}^1(\mathcal{U})$ is defined as
$$
\{(f,\omega)\in C^{1}(\mathcal{U},\mathcal{O}_Y)\times C^0(\mathcal{U},\Omega_Y^1)\, |\, f_{i,k}=f_{i,j}+f_{j,k},\diff f_{i,j}=\omega_i-\omega_j, \diff \omega_i=0 \}
$$ 
and $B_{dR}^1(\mathcal{U})$ is defined as
$$
\{(f,\omega)\in C^{1}(\mathcal{U},\mathcal{O}_Y)\times C^0(\mathcal{U},\Omega_Y^1)\, |\,\exists g\in C^0(\mathcal{U},\mathcal{O}_Y),f_{i,j}=g_i-g_j,  \omega_i=\diff g_i \}\, .
$$
The first de Rham cohomology $H^1_{dR}(Y)$ is independent of the choice of the open affine cover and we have 
$$
H^1_{dR}(Y)\cong H^1_{dR}(\mathcal{U})\, .
$$
 In particular for the case of a curve $X$,  we can take an open affine cover $\mathcal{U}=\{U_1,U_2\}$ consisting of only two open parts. Then the de Rham cohomology $H^1_{dR}(Y)$ can be described using $(\ref{de Rham defi})$
with
$$
Z_{dR}^1(\mathcal{U})=\{(t,\omega_1,\omega_2)|~t\in \mathcal{O}_Y(U_1\cap U_2),\omega_i\in \Omega_Y^1(U_i),\diff t=\omega_1-\omega_2\}
$$
and 
$$
B_{dR}^1(\mathcal{U})=\{(t_1-t_2,\diff t_1,\diff t_2)|~t_i\in \mathcal{O}_Y(U_i)\}\, .
$$
\subsection{The Verschiebung operator and the Cartier operator}
The Cartier operator $\mathcal{C}$ is the operator on the rational differential forms on a curve $X$ defined in \cite{MR0084497}.  If $x$ is a separating variable of $k(X)$, any $f\in k(X)$ can be written as 
$$
f=f_0^p+f_1^px+\dots+f_{p-1}^px^{p-1}, \, \text{with} f_i\in k(X)\, ,
$$
then for a rational differential form $\omega=f\diff x$ with $f$ as above, we have $\mathcal{C}(\omega)=f_{p-1}\diff x$. 

We define the \textit{rank of the Cartier operator} to be the rank of the Cartier operator restricted to $H^0(X,\Omega_X^1)$ and the \textit{Cartier-Manin matrix} is the matrix of the Cartier operator with respect to a given basis.

Suppose $A$ is a principally polarized abeian variety of dimension $g$ defined over $k$, e.g. $A={\rm Jac}(X)$ for some curve $X$ of genus $g$. Let $F_{A/k}:A\mapsto A^{(p)}$ be the relative Frobenius morphism of $A$ defined by the universal property of  $A^{(p)}$, which is the pullback of $A\mapsto {\rm Spec}(k)$ via the absolute Frobenius map  $\sigma: {\rm Spec}(k) \mapsto {\rm Spec}(k)$. The relative Frobenius morphism of $A$ is purely inseparable of degree $p^g$. Note that the multiplication by $p$ on $A$ is of degree $p^2g$ and it factors through the relative Frobenius morphism, i.e. $[p]_A=V_{A/k}\circ F_{A/k}$ with $V_{A/k}: A^{(p)}\mapsto A$ the Verschiebung morphism. 

Note that  $V_{A/k}$ induces a  $\sigma^{-1}$-linear map of $k$-vector spaces 
$$
V: H^1_{dR}(A)\mapsto H^1_{dR}(A)
$$  
by considering the fact $H^1_{dR}(A^{(p)})\cong(k,\sigma)\otimes_k H^1_{dR}(A)$. Moreover, for a curve $X$ and $(f,\omega)\in H^1_{dR}(X)$, one has $V(f,\omega)=(0,V(\omega))=(0,\mathcal{C}(\omega))$. Similarly one have the Frobenius operator $F:H^1_{dR}(A)\mapsto H^1_{dR}(A)$ induced by the relative Frobenius of $A$.

\subsection{The Ekedahl-Oort stratification}\label{E-O subsection}
The Ekedahl-Oort stratification is defined by a final filtration on the first de Rham cohomology of abelian varieties, which we define as follows.
\begin{defn}
Let $(A,\lambda)$ be a principally polarized abelian variety with $G=H^1_{dR}(A)$. A \textit{final filtration} is a filtration 
\begin{align*}
0= G_0\subset G_1\subset \dots\subset G_g=V(G)\subset \dots\subset G_{2g}=G
\end{align*}
which is stable under $V$ and $F^{-1}$ with the following properties being satisfied:
\begin{align*}
\dim (G_i)=i,\, G_{i}^{\perp}=G_{2g-i}, i=1,\dots,2g.
\end{align*}
The associated \textit{final type} $v$ is the increasing and surjective map 
\begin{align*}
v:\{0,1,2,\dots,2g\}\mapsto \{0,1,\dots,g\}
\end{align*}
such that $ V(G_i)=G_{v(i)}$.
\end{defn}
By the properties of $V$ and $F$, we have $v(2g-i)=v(i)-i+g$ for $0\leq i\leq g$.
To find a final filtration of $(A,\lambda)$, we first construct a weak version of it, namely the \textit{canonical filtration}, and refine it to a final one. 

Now to construct the canonical filtration of $A$, we start with $0\subset G$ and add the images of $V$ 
and the orthogonal complements of the images.  Repeating this process for a finite number of  times, we obtain a filtration
$$
0= C_0\subset G_1\subset \dots\subset G_r=V(G)\subset \dots\subset G_{2r}=G\, ,
$$
 which is stable under $V$ and $\, ^\perp$. Then it is also stable under $F^{-1}$ and it is unique. This filtration is called canonical filtration and can be refined to a final one by  choosing a filtration of length $2g$ which is stable under $V$ and $\, ^\perp$. In general a final filtration is not unique but the final type is unique, cf. \cite{Oort1999}. 
 
Given a final type $v$, we associate to it a Young diagram, or  equivalently  a $n$-tuple
$\mu=[\mu_1, \ldots,\mu_n]$ with $0\leq n \leq g$ and $\mu_1 > \mu_2
>\cdots >\mu_n>0$ such that 
$$
\mu_j=\#\{i:\, 1\leq i\leq j,v(i)+j\leq i   \}\, .
$$
Denote by $\mathcal{Z}_{\mu}$ the set of  geometric points of $\mathcal{A}_g$ with given final type $v$ with Young diagram $\mu$. We say an abelian variety (resp. a curve) has Ekedahl-Oort type $\mu$ if the corresponding point in $\mathcal{A}_g$ (resp. $\mathcal{M}_g$) lies in $\mathcal{Z}_{\mu}$ (resp. $Z_{\mu}$).
\subsection{The $p$-rank and $a$-number}
Let $X$ be a  curve of genus $g$ defined over $k$.
 Such curve has several invariants, e.g. the $a$-number and the $p$-rank.
The $a$-number of the curve $X$ is defined as $a_X=\dim_k({\rm Hom}(\alpha_p,{\rm Jac}(X)))$ with $\alpha_p$ the group scheme which is the  kernel of Frobenius on the additive group scheme $\mathbb{G}_a$.  The $a$-number of $X$ is equal to $g-r$ where $g$ is the genus of $X$ and  $r$ is the rank of the Cartier-Manin matrix, cf. \cite{MR0084497,Seshadri1958-1959}. The $p$-rank of a curve $X$ is the number $f_X$ such that $\#{\rm Jac}(X)[p](k)=p^{f_X}$. One see that $1\leq a_X+f_X\leq g$. Moreover, if $X$ has Ekedahl-Oort type $\mu=[\mu_1,\dots,\mu_n]$, the $f_X=g-\mu_1$ and $a_X=n$.

\section{Proof of Theorem \ref{thmmain}} 
\subsection{Set up}Recall that  we denote by $\mathcal{M}_g$   the
moduli space (stack) of (smooth projective) curves of genus $g$ defined over $k$ and by $\overbar{\mathcal{M}}_g$ its Deligne-Mumford compactification.
Write $\overbar{\mathcal{Z}_{\mu}}$ for the Zariski closure of $\mathcal{Z}_{\mu}$ in $\mathcal{A}_g$. 
Put $\Wtilde_g:=\overbar{\mathcal{M}}_g-\Delta_0$. It is  an open substack of $\overbar{\mathcal{M}}_g$ parametrising stable curves of compact type of genus $g$. Note that a stable curve is said to be of compact type if the dual graph of the curve is a tree, we refer the definition of the dual graph of a curve to \cite{MR1827016}. For $g=1$, a stable curve of genus $1$ of compact type is equivalent to an irreducible smooth elliptic curve.   For $g\geq 2$, the Jacobian of  a stable curve of genus $g\geq 2$ with a tree as its dual graph  is an abelian variety.

The Torelli morphism ${\tau}:{\mathcal{M}}_g\mapsto \mathcal{A}_g$ can be extended to a regular map:
$$
\tilde{\tau}:\, \Wtilde_g\, \mapsto \mathcal{A}_g	\, .
$$
Since curve with different marked point can have  same Jacobian, the map $\tilde{\tau}$ has positive dimensional fibres.
It is known that this morphism is proper and  its image $\tilde{\tau}(\Wtilde_g)$ is a reduced closed subscheme of $\mathcal{A}_g$, cf. \cite{MR3184184}.  For $g=4$, we have that $\dim (\mathcal{M}_4)=9$ and $\dim(\mathcal{A}_4)=10$. Moreover, by a work of Igusa \cite{MR656035}, the image $\tilde{\tau}(\Wtilde_4)$ is an ample divisor on $\mathcal{A}_4$. Hence $\tilde{\tau}(\Wtilde_4)\cap \mathcal{Z}_{\mu}$ is of codimension at most $\sum_{i=1}^n \mu_i$ with $\mu=[\mu_1,\dots,\mu_n]$. Recall that there is a partial order $\preceq$ on the set of Ekedahl-Oort strata. We have $\overbar{Z_{\mu}}\supset \overbar{Z_{\upsilon}}$ if $\mu\preceq \upsilon$. 

Now we can give a proof of Theorem $\ref{thmmain}$.

\begin{proof}[Proof of Theorem $\ref{thmmain}$]
The $p$-rank part were known by \cite[Theorem 2.3]{MR2084584}, saying that $Z_{\mu}$ is of codimention $\mu_1$ in $\mathcal{M}_g$ for $\mu=[\mu_1]$.
For $\mu=[2,1]$, we compute the dimension of $\tilde{\tau}(\Delta_i)\cap \mathcal{Z}_{[2,1]}$ in ${\mathcal{A}_4}$ for $i=1,2$.
A curve corresponding to a point in the pull back of a generic point in $\tilde{\tau}(\Delta_1)\cap \mathcal{Z}_{[2,1]}$ is formed from two smooth pointed curves $X_1$ and $X_3$ of genus 1 and 3 respectively. Then $X_1$ is ordinary and $X_3$ has Ekedahl-Oort type $\mu=[2,1]$ or $X_1$ is supersingular and $X_3$ has Ekedahl-Oort type $\mu=[1]$. In the first case, we consider the action of Verschiebung operator $V$ on the de Rham cohomology. Let $A_2=H^1_{dR}(X_1), B_6=H^1_{dR}(X_3)$ and $C_8=H^1_{dR}(X_1\times X_3)$. We now compute the final filtration on $H^1_{dR}(X_1\times X_3)$ as explained in Section $\ref{E-O subsection}$. Recall that we have a pairing $\langle\, ,\, \rangle$ on $H^1_{dR}(X)$ for any smooth irreducible projective curve $X$. We have
$$
V(A_2)=A_1,\,V(A_1)=A_1,\,A_1^{\perp}=A_1\, ,
$$
where $A_1^{\perp}$ is the orthogonal complement of $A_1$ in $A$ with respect to the pairing $\langle\, ,\, \rangle$.
Similarly, 
$
V(B_6)=B_3,V(B_3)=B_1,V(B_1)=B_1,\, B_1^{\perp}=B_5
$.
For a final filtration $
0\subset G_1\subset G_2\subset \dots \subset G_g\subset \dots \subset  G_{2g}
$, we have 
\begin{align}
\dim V(G_{2g-n})=\dim V(G_{n})+g-n \label{property of canonical type}
\end{align}
with $0\leq n\leq g$. Then $\dim V(B_5)=\dim V(B_1)+3-1=3$. It follows that $V(B_5)=B_3$ and
\begin{align*}
V(C_8)&=V(A_2+B_6)=A_1+B_3=C_4,\,\\ V(C_4)&=V(A_1+B_3)=A_1+B_1=C_2,\\\, V(C_2)&=V(A_1+B_1)=A_1+B_1\, .
\end{align*}
Furthermore,
$$
C_2^{\perp}=A_1+B_5=C_6, V(C_6)=V(A_1+B_5)=A_1+B_3=C_4\, .
$$
This gives the canonical filtration $
0\subset C_2\subset C_4\subset C_6\subset C_8
$ for $C_8=H^1_{dR}(X_1\times X_3)$ and the action of $V$.
So the Ekedahl-Oort type is $[2,1]$. Note that the moduli space of pointed elliptic curves  is irreducible of dimension 1 and $Z_{[2,1]}\subset \mathcal{M}_3$ is irreducible of dimension 3 by \cite[Theorem 11.3]{Ekedahl2009}. Hence in this case the component in $\tilde{\tau}(\Delta_1)\cap \mathcal{Z}_{[2,1]}$ is of dimension 4. By a similar argument, in the other case, a generic point also has  Ekedahl-Oort type $[2,1]$ and the component of $\tilde{\tau}(\Delta_1)\cap \mathcal{Z}_{[2,1]}$ is non-empty of dimension 5. Hence $\tilde{\tau}(\Delta_1)\cap \mathcal{Z}_{[2,1]}$ is of dimension at most 5.

 For $\tilde{\tau}(\Delta_2)\cap \mathcal{Z}_{[2,1]}$,  a generic point in $\tilde{\tau}(\Delta_2)\cap \mathcal{Z}_{[2,1]}$ is the union of two pointed curves $X_2$ and $\check{X}_2$, both of genus 2, where $X_2$ and $\check{X}_2$ both have $3$-rank $1$ or $X_2$ is ordinary and $\check{X}_2$ is superspecial. In the first case, let $A_4=H^1_{dR}(X_2), B_4=H^1_{dR}(\check{X}_2)$ and $C_8=H^1_{dR}(X_2\times \check{X}_2)$, then we have
$$
V(A_4)=A_2,\,V(A_2)=A_1,\,A_1^{\perp}=A_3\, . 
$$
By equation $(\ref{property of canonical type})$, we have $V(A_3)=A_2$. We also have $V$ acts similarly on $B_4$.
Furthermore, 
\begin{align*}
V(C_8)=V(A_4+B_4)=A_2+B_2=C_4,\\
V(C_4)=V(A_2+B_2)=A_1+B_1=C_2,\\
V(C_2)=V(A_1+B_1)=A_1+B_1\, .
\end{align*}
Combined with $C_2^{\perp}=A_3+B_3=C_6$ and $V(C_6)=A_2+B_2=C_4$, we have the canonical filtration $
0\subset C_2\subset C_4\subset C_6\subset C_8$ and the action of $V$.
 Hence the generic point has Ekedahl-Oort type $[2,1]$ and the component is of dimension 4. In a similar fashion we can prove in the other case that the component is of dimension 3.

 On the other hand, $\tilde{\tau}(\Wtilde_4)$ is ample in $\mathcal{A}_4$ and hence $\tilde{\tau}(\Wtilde_4)\cap \mathcal{Z}_{[2,1]}$ is of dimension at least $6$. Furthermore, $\mathcal{Z}_{[2,1]}\subset \mathcal{A}_4$ is irreducible of dimension 7 by \cite[Theorem 11.5]{Ekedahl2009}. Suppose $Z_{[2,1]}\subset \mathcal{M}_4$ is of dimension 7. By taking the closure in $\mathcal{A}_4$ and by the irreducibility of $\mathcal{Z}_{[2,1]}$,  we have $\overbar{\mathcal{Z}_{[2,1]}}=\overbar{\tilde{\tau}(Z_{[2,1]})}$. Since $\tilde{\tau}(\Wtilde_4)$ is closed in $\mathcal{A}_4$, we have  $\overbar{\mathcal{Z}_{[2,1]}}\subset \tilde{\tau}(\Wtilde_4)$. This implies  $\overbar{\mathcal{Z}_{[4,1]}}\subset \overbar{\mathcal{Z}_{[2,1]}}\subset \tilde{\tau}(\Wtilde_4)$ with $\overbar{\mathcal{Z}_{[4,1]}}$  irreducible of dimension $5$ by \cite[Theorem 11.3]{Ekedahl2009}. Note that we also have $\overbar{\mathcal{Z}_{[4,1]}}\subset \overbar{\mathcal{Z}_{[4]}}$. Then $\overbar{\mathcal{Z}_{[4,1]}}\subset \overbar{\mathcal{Z}_{[4]}}\cap \tilde{\tau}(\Wtilde_4)$ with $\dim(\overbar{\mathcal{Z}_{[4,1]}})=\dim(\overbar{\mathcal{Z}_{[4]}}\cap \tilde{\tau}(\Wtilde_4))=5$ by Theorem \ref{thm V_f codimension}. Then $\overbar{\mathcal{Z}_{[4,1]}}$ is an irreducible component of $\overbar{\mathcal{Z}_{[4]}}\cap \tilde{\tau}(\Wtilde_4)$ and every abelian variety corresponds to a point of this component of $\overbar{\mathcal{Z}_{[4]}}\cap \tilde{\tau}(\Wtilde_4)$ has $a$-number $\geq 2$. Write $V_0(\Wtilde_4)$ for the locus of curves with $p$-rank 0 in $\Wtilde_4$. By the proof  of Theorem 2.3 in  \cite{MR2084584},  a generic point of any irreducible component of the locus  $\tilde{\tau}(V_0(\Wtilde_4))$  has $a$-number 1. This contradiction implies  $\overbar{\mathcal{Z}_{[2,1]}}\not\subset \tilde{\tau}(\Wtilde_4)$ and $Z_{[2,1]}$ is non-empty of dimension $6$.

 For $\mu=[3,1]$, similarly we can compute that $\tilde{\tau}(\Delta_i)\cap \mathcal{Z}_{[3,1]}$ has dimension at most 3 for $i=1,2$.  On the other hand, $\tilde{\tau}(\Wtilde_4)\cap \mathcal{Z}_{[3,1]}$ is of dimension at least 5  as   $\tilde{\tau}(\Wtilde_4)$ in $\mathcal{A}_4$ is ample. Moreover,  $\mathcal{Z}_{[3,1]}\subset \mathcal{A}_4$ is irreducible of dimension 7 by \cite[Theorem 11.5]{Ekedahl2009}.  Suppose $Z_{[3,1]}$ is of dimension $6$. By taking the closure in $\mathcal{A}_4$ and by the irreducibility of $\mathcal{Z}_{[3,1]}$,  we have $\overbar{\mathcal{Z}_{[3,1]}}=\overbar{\tilde{\tau}(Z_{[3,1]})}$. Note that $\overbar{\mathcal{Z}_{[4,1]}}\subset\overbar{\mathcal{Z}_{[3,1]}}$. By the similar argument in the case $Z_{[2,1]}$, we have  $\overbar{\mathcal{Z}_{[3,1]}}\not\subset \tilde{\tau}(\Wtilde_4)$ and the locus $Z_{[3,1]}$ is non-empty of dimension 5.

For $\mu=[3,2]$, we show that a smooth curve with affine equation 
\begin{align}
y^3+y^2+by=x^5+a_3x^3+a_2x^2+a_0\, ,\label{equation E-O [3,2]}
\end{align}
where $a_i,b\in k$ and $a_2\neq 0$  has Ekedahl-Oort type $[3,2]$. 

We first show that the map from the parameter space with coordinates $(a_3,a_2,a_0,$ $b)$ to $\mathcal{M}_4$ has finite fibres. Denote by $\sigma$ an isomorphism between two smooth Artin-Schreier curves given by $y^3+y^2+by=x^5+a_3x^3+a_2x^2+a_0$ and $y^3+y^2+b_1y=x^5+c_3x^3+c_2x^2+c_0$ as in $(\ref{example [3,2,1]})$. After possibly composing with an inversion $x\mapsto 1/x$, we may assume $\sigma(x)=\alpha x+\beta$. Also since $\sigma$ is invertible, we have $\sigma (y)=zy+\delta$ with $z$ a unit in $k$ and $\delta\in k(x)$. Hence
$$
z^3y^3+\delta^3+(zy+\delta)^2+b(zy+\delta)=(\alpha x+\beta)^5+a_3(\alpha x+\beta)^3+a_2(\alpha x+\beta)^2+a_0\, .
$$
Then we have $\delta=\beta=0$ and $y^3+y^2/z+by/z^2=1/z^3(\alpha^5x^5+a_3\alpha^3x^3+a_2\alpha^2x^2+a_0)$. This implies $\beta=\delta=0,\alpha^5=z=1$.

Now we show that a curve $X$ with equation $(\ref{equation E-O [3,2]})$ has Ekedahl-Oort type $[3,2]$. A basis of 
$H^0(X,\Omega_X^1)$ is given by $\omega_1=1/(y-b)\diff x, \omega_2=x/(y-b) \diff x, \omega_3=x^2/(y-b)\diff x,\omega_4=1\diff x$. Then the Cartier-Manin matrix is given by
$$
\left ( \begin{array}{lccl}
a_2 & 0 & a_0-b^3 &0 \\ 
1 & 0 & a_3 &0\\
0 & 0 & 0 &0\\
0 & 0 & 1 &0\\
\end{array}\right )^{1/3}\, 
$$
and it has  rank $2$ and semisimple rank $1$. One can compute that curve $X$ has final filtration 
$$
0\subset G_1\subset G_2\subset G_3\subset G_4=H^0(X,\Omega_X^1)\subset \dots \subset G_8=H^1_{dR}(X)\, ,
$$
where $G_1=\langle a_2^{1/3}\omega_1+\omega_2\rangle, G_2=\langle a_2\omega_1+\omega_2, (a_0-b^3)^{1/3}\omega_1+a_3^{1/3}\omega_2+\omega_4\rangle$ and $G_3=\langle \omega_1,\omega_2,\omega_4\rangle$. Hence the curve has Ekedahl-Oort type $[3,2]$ and $ Z_{[3,2]}$ is non-empty of dimension at least 4. 

%Note that again $\mathcal{Z}_{[3,2]}$ is irreducible in $\mathcal{A}_4$ by \cite[Theorem 11.5]{Ekedahl2009}. By a similar argument in case $Z_{[2,1]}$, we have $\overbar{\mathcal{Z}_{[3,2]}}\not\subset \tilde{\tau}(\Wtilde_4)$.

For $\mu=[3,2,1]$, we show that the family of curves given by equation 
\begin{align}
y^3-bx^3(y^2+y)=x^5+cx^3+dx^2+1, ~b,c,d \in k,~bd\neq 0 \label{example [3,2,1]}
\end{align}
has Ekedahl-Oort type $[3,2,1]$. 

In a similar fashion, we first show that the map from the parameter space with coordinates $(b,c,d)$ to $\mathcal{M}_4$ has finite fibres. Denote by $\sigma$ an isomorphism between two smooth Artin-Schreier curves given by $y^3-bx^3(y^2+y)=x^5+cx^3+dx^2+1$ and $y^3-b_1x^3(y^2+y)=x^5+c_1x^3+d_1x^2+1$ as in $(\ref{example [3,2,1]})$. After possibly composing with~an~inversion $x\mapsto 1/x$, we may assume $\sigma(x)=\alpha x+\beta$. Also since $\sigma$ is invertible, we have $\sigma (y)=zy+\delta$ with $z$ a unit in $k$ and $\delta\in k(x)$. Hence
$$
z^3y^3+\delta^3+b(\alpha x+\beta)^3(z^2y^2+(2\delta+z)y+\delta^2+\delta)=\sigma(x)^5+c\sigma(x)^3+d\sigma(x)^2+1\, .
$$
By comparing the coefficients of $x^i$ for $i=0,1,2,3,4,5$, we have $\beta=\delta=0$ and $z=1$. Furthermore, we have $y^3-b\alpha^3x^3(y^2+y)=\alpha^5x^5+\alpha^2 x^3+dx^2+1$ and hence $\alpha^5=1$. This implies $\beta=\delta=0,\alpha^5=z=1$.

Now we compute the Ekedahl-Oort type of a curve $X$ given by equation $(\ref{example [3,2,1]})$. For a basis of $H^0(X,\Omega_X^1)$ we choose
$$1/s(x,y)\diff x, x/s(x,y)\diff x, x^2/s(x,y)\diff x, y/s(x,y) \diff x$$
with $s(x,y)=x^3(y-1)$.
Then the Cartier-Manin matrix has rank 1 and semisimple rank $1$.
Hence $ Z_{[3,2,1]}$ is non-empty with dimension at least 3.

% Similar to an argument in case $\mu=[3,2]$,  the locus $Z_{[3,2,1]}$ is irreducible \cite[Theorem 11.5]{Ekedahl2009} in $\mathcal{M}_4$ and $\overbar{\mathcal{Z}_{[3,2,1]}}\not\subset \tilde{\tau}(\Wtilde_4)$. 

For $\mu=[4,1]$, we have $\mathcal{Z}_{[4,1]}$ is irreducible by \cite[Theorem 11.5]{Ekedahl2009}.  For $\tilde{\tau}(\Delta_1)\cap \mathcal{Z}_{[4,1]}$, a curve corresponding to a point in the pull back of a generic point in $\tilde{\tau}(\Delta_1)\cap \mathcal{Z}_{[4,1]}$ is formed from two pointed curves $X_1$ and $X_3$ of genus 1 and 3 respectively. Moreover, $X_1$ is superspecial and $X_3$ has Ekedahl-Oort type $\mu=[3]$. One can easily compute that in this case the curve formed from pointed curves $X_1$ and $X_3$ has Ekedahl-Oort type $[4,2]$ and  $\tilde{\tau}(\Delta_1)\cap \mathcal{Z}_{[4,1]}$ is of dimension $3$.  For $\tilde{\tau}(\Delta_2)\cap \mathcal{Z}_{[4,1]}$, a point in the pull back of a generic point in $\tilde{\tau}(\Delta_2)\cap \mathcal{Z}_{[4,1]}$ is formed from two pointed  curves with Ekedahl-Oort type $\mu=[2]$. Then one can compute that a  generic point of $\tilde{\tau}(\Delta_2)\cap \mathcal{Z}_{[4,1]}$ has Ekedahl-Oort type $[4,3]$ and $\tilde{\tau}(\Delta_2)\cap \mathcal{Z}_{[4,1]}$ is of dimension $2$. 

On the other hand, $\tilde{\tau}(\Wtilde_4)\cap \mathcal{Z}_{[4,1]}$ is of dimension at least $4$. Hence $Z_{[4,1]}$ is of dimension at least 4. Note that if $Z_{[4,1]}$ is of dimension 5. Then by taking the closure in $\mathcal{A}_4$ and by the irreducibility of $\mathcal{Z}_{[4,1]}$, we have $\overbar{\mathcal{Z}_{[4,1]}}\subset \tilde{\tau}(\Wtilde_4)$, which is a contradiction by the proof of the case $Z_{[2,1]}$.
 Hence $Z_{[4,1]}$ is non-empty of dimension 4.

For $\mu=[4,2]$ and $[4,3]$, we have $\mathcal{H}_4\cap Z_{\mu}$ non-empty hence $ Z_{\mu}$ is non-empty by  \cite[Theorem 1.1]{2018arXiv181204996Z}.

Denote $\theta(x)=x^3+x^2+b_1x+b_2$ and $h(x)=x^5+a_4x^4+a_3x^3+a_2x^2+a_1x$. We consider the family of curves parametrized by equation
\begin{align}
y^3+y^2+\theta(x)y=h(x)\, ,
 \label{example 2-dim [4,3] 3-dim}
\end{align}
where $a_i,b_j\in k$ such that $b_2=b_1^2+a_2,~a_4=b_1-1$ and $a_3=b_1^3+b_1^2+b_1+a_2$. 

Note that the map from the parameter space to the $\mathcal{M}_4$ is finite and hence it gives a $3$-dimensional sublocus in $\mathcal{M}_4$. Indeed,
denote by $\sigma$ an isomorphism between two smooth Artin-Schreier curves given by $\theta_1(x),h_1(x)$ and $\theta_2(x),h_2(x)$ as in $(\ref{example 2-dim [4,3] 3-dim})$. After possibly composing with an inversion $x\mapsto 1/x$, we may assume $\sigma(x)=\alpha x+\beta$. Also since $\sigma$ is invertible, we have $\sigma (y)=zy+\delta$ with $z$ a unit in $k$ and $\delta\in k(x)$. Then one can easily show that $z=\alpha=1,~\beta^3+\beta^2=0$ and $-(\delta^3+\delta^2+b_2\delta)+h_1(\beta)=0$ with $\theta_1(x)=x^3+x^2+b_1x+b_2$ and $h(x)=x^5+a_4x^4+a_3x^3+a_2x^2+a_1x$.

Indeed, denote by $\sigma$ an isomorphism between two smooth Artin-Schreier curves given by $y^3-y=f_1(x)=\sum_{i=0}a_ix^i$ and $y^3-y=f_2(x)=\sum_{i=0}b_ix^i$ as in $(\ref{example 2-dim [4,3]})$. By \cite[Lemma $2.1.5$]{10.2307/25099135}, two Artin-Schreier curves are isomorphic if and only if $f_2(x)=z f_1(x)+\delta^3-\delta$ with $z\in \mathbb{F}_3^*$ and $\delta\in k(x)$. Moreover, the $\sigma$ is defined by $\sigma(x)=\alpha
 x+\beta$ and $\sigma (y)=zy+\delta$. Hence we have $\beta=0, \alpha^5=z\in \mathbb{F}_3^*$ and $\delta\in \mathbb{F}_3$.

Now we show that any smooth curve $X$ given by equation $(\ref{example 2-dim [4,3] 3-dim})$ above has Ekedahl-Oort type $[4,3]$.
For a basis of $H^0(X,\Omega_X^1)$, we choose 
$$1/s(x,y)\diff x, x/s(x,y)\diff x, x^2/s(x,y)\diff x, y/s(x,y) \diff x$$
with $s(x,y)=y-(x^3+x^2+b_1x+b_2)$.  For the Cartier operator $\mathcal{C}$ we have 
\begin{align*}
\mathcal{C}(\frac{1}{s(x,y)}\diff x)&=\mathcal{C}(\frac{1}{y-\theta(x)})=\mathcal{C}(\frac{(y-\theta (x))^2}{(y-\theta (x))^3})\\
&=\frac{1}{y-\theta (x)}\mathcal{C}(y^2+\theta ^2(x)+\theta(x) y)\, .
\end{align*}
Note that $y^2+\theta (x)y=-y^3+h(x)$ and $\theta ^2(x)=x^6+2x^5+(2b_1+1)x^4+(2b_1+2b_2)x^3+(b_1^2+2b_2)x^2+2b_1b_2x+b_2^2$.
Hence 
\begin{align*}
y^2+\theta ^2(x)+\theta(x) y=-y^3+x^6+(2b_1+1+a_4)x^5+(2b_1+2b_2+a_3)x^3\\
+(b_1^2+2b_2+a_2)x^2+(2b_1b_2+a_1)x+b_2^2\, ,
\end{align*}
and we have \begin{align*}
\mathcal{C}(\frac{1}{s(x,y)}\diff x)&=\frac{1}{y-\theta (x)}\mathcal{C}(y^2+\theta ^2(x)+\theta(x) y)\\
&=\frac{(2b_2+b_1^2+a_2)^{1/3}}{y-\theta(x)}=(2b_2+b_1^2+a_2)^{1/3}\frac{1}{s(x,y)}\diff x\, .
\end{align*}
Similarly one can compute $\mathcal{C}(x/s(x,y)),\mathcal{C}(x^2/s(x,y))$ and $\mathcal{C}(y/s(x,y))$. Then we have the Cartier-Manin matrix which equals to 
\begin{align*}
&\left ( \begin{array}{cccc}
2b_2+b_1^2+a_2 & 2b_1b_2+a_1 & b_2^2 & b_2^2+b_1(2b_1b_2+a_1)\\ 
0 & 2b_1+1+a_4 & 2b_1+2b_2+a_3 & 2b_1+2b_2+a_3\\
0 & 0 & 1 &1\\
0 & 0 & -1 &-1\\
\end{array}\right )^{1/3}\\
&=\left ( \begin{array}{cccc}
0 & 2b_1b_2+a_1 & b_2^2 & b_2^2+b_1(2b_1b_2+a_1)\\ 
0 & 0 & b_1^3 & b_1^3\\
0 & 0 & 1 &1\\
0 & 0 & -1 &-1\\
\end{array}\right )^{1/3}\, .
\end{align*}
We obtain that ${\rm rank}(\mathcal{C})=2$ and $\mathcal{C}^2=0$ on $H^0(X,\Omega_X^1)$.
Hence $X$ with equation $(\ref{example 2-dim [4,3] 3-dim})$ has Ekedahl-Oort type $[4,3]$ and $Z_{[4,3]}$ is of dimension $3$.

For $\mu=[4,3,2,1]$, a curve corresponding to a point in $\mathcal{Z}_{[4,3,2,1]}$ is superspecial.  Then $Z_{[4,3,2,1]}$ is empty as Ekedahl   \cite{Ekedahl1987} showed that any superspecial curve of genus $g$ satisfies :
$$
g\leq 3(3-1)/2=3\, .
$$

\end{proof}
\begin{rmk}
As an example,  a family of curves given by the equation 
$$
y^3+bx^2y=x^5+ax^4+x,~a,b\in k, b\neq 0
$$ 
has Ekedahl-Oort type $[2,1]$.
Moreover, the family of curves parametrized by equation
\begin{align}
y^3-y=x^5+a_2x^2+a_1x,~a_1,a_2\in k\,  \label{example 2-dim [4,3]}
\end{align}
gives a $2$-dimensional locus in $ Z_{[4,3]}$. 
 On the other hand, the family of curves parametrized by equation
\begin{align}
y^3+y^2+(x^3+x^2)y=x^5+2x^4+a_1x+a_2\, ,
 \label{example 2-dim [4,3] second}
\end{align}
where $a_1,a_2\in k$ and $a_1\neq 0$, gives another $2$-dimensional locus in $ Z_{[4,3]}$. Consider a curve with equation $(\ref{example 2-dim [4,3] second})$ above, with respect to the basis given by 
$$1/s(x,y)\diff x, x/s(x,y)\diff x, x^2/s(x,y)\diff x, y/s(x,y) \diff x$$
with $s(x,y)=y+2x^3+2x^2$, one can easily compute the Cartier-Manin matrix which equals to
$$
\left ( \begin{array}{lccl}
0 & a_1 & a_2 & a_2\\ 
0 & 0 & 0 &0\\
0 & 0 & 1 &1\\
0 & 0 & -1 &-1\\
\end{array}\right )^{1/3}\, .
$$
Then the rank and semisimple rank of the Cartier-Manin matrix is 2 and 0 respectively and the curve with equation $(\ref{example 2-dim [4,3] second})$ has Ekedahl-Oort type $[4,3]$.

\end{rmk}
\begin{rmk}
Note that Achter and Pries \cite{MR3287680} have some result about generic Newton polygon of curves of given genus and $p$-rank.
\end{rmk}
\section*{Acknowledgement}{I would like to thank my supervisor, Professor Gerard van der Geer,  for his patient and continuous guidance and kind advices in every stage of this paper. I also would like to thank Katsura and Pries for insightful comments.}

%\bibliography{/Users/zhouzijian/Dropbox/bibfile}

\begin{thebibliography}{10}
\bibitem{MR3287680}
J.~Achter and R.~Pries.
\newblock {\em Generic {N}ewton polygons for curves of given {$p$}-rank.
  \upshape {I}n: Algebraic curves and finite fields}, pages 1--21.
\newblock De Gruyter, Berlin, 2014.

\bibitem{MR1827016}
F.~Andreatta.
\newblock {\em On {M}umford's uniformization and {N}\'eron models of
  {J}acobians of semistable curves over complete rings. \upshape In: {Moduli of
  abelian varieties}}, pages 11--126.
\newblock Birkh\"auser, Basel, 2001.

\bibitem{MR0084497}
P.~Cartier.
\newblock Une nouvelle op\'eration sur les formes diff\'erentielles.
\newblock {\em C. R. Acad. Sci. Paris}, 244:426--428, 1957.

\bibitem{2017arXiv171204921D}
S.~{Devalapurkar} and J.~{Halliday}.
\newblock {The Dieudonn\'e modules and Ekedahl-Oort types of Jacobians of
  hyperelliptic curves in odd characteristic}.
arXiv:1712.04921, Dec 2017.

\bibitem{EGA}
J.~Dieudonn{\'e} and A.~Grothendieck.
\newblock \'{E}l\'ements de g\'eom\'etrie alg\'ebrique.
\newblock {\em Inst. Hautes \'Etudes Sci. Publ. Math.}, 1961--1967.

\bibitem{Ekedahl1987}
T.~{Ekedahl}.
\newblock On supersingular curves and abelian varieties.
\newblock {\em Math. Scand.}, 60:151--178, 1987.

\bibitem{Ekedahl2009}
T.~Ekedahl and G.~van~der Geer.
\newblock {\em Cycle Classes of the E-O Stratification on the Moduli of Abelian
  Varieties. \upshape {I}n: {Algebra, Arithmetic, and Geometry: Volume I: In
  Honor of Yu. I. Manin}}, pages 567--636.
\newblock Birkh{\"a}user Boston, Boston, 2009.

\bibitem{MR3095219}
A.~Elkin and R.~Pries.
\newblock Ekedahl-{O}ort strata of hyperelliptic curves in characteristic 2.
\newblock {\em Algebra Number Theory}, 7(3):507--532, 2013.

\bibitem{MR2084584}
C.~Faber and G.~van~der Geer.
\newblock Complete subvarieties of moduli spaces and the {P}rym map.
\newblock {\em J. Reine Angew. Math.}, 573:117--137, 2004.

\bibitem{MR0463157}
R.~Hartshorne.
\newblock {\em Algebraic geometry}, volume~52 of {\em Graduate Texts in
  Mathematics}.
\newblock Springer-Verlag, New York-Heidelberg, 1977.

\bibitem{MR656035}
J.~Igusa.
\newblock On the irreducibility of {S}chottky's divisor.
\newblock {\em J. Fac. Sci. Univ. Tokyo Sect. IA Math.}, 28(3):531--545 (1982),
  1981.

\bibitem{MR3184184}
B.~Moonen and F.~{Oort}.
\newblock {\em The {T}orelli locus and special subvarieties. \upshape {I}n:
  Handbook of moduli. {V}ol. {II}}, pages 549--594.
\newblock Int. Press, Somerville, MA, 2013.

\bibitem{MR0241435}
T.~Oda.
\newblock The first de {R}ham cohomology group and {D}ieudonn\'e modules.
\newblock {\em Ann. Sci. \'Ecole Norm. Sup. (4)}, 2:63--135, 1969.

\bibitem{Oort1999}
F.~Oort.
\newblock {\em A Stratification of a Moduli Space of Polarized Abelian
  Varieties in Positive Characteristic. \upshape In: {Moduli of Curves and
  Abelian Varieties: The Dutch Intercity} {Seminar on Moduli}}, pages 47--64.
\newblock Vieweg+Teubner Verlag, Wiesbaden, 1999.

\bibitem{10.2307/25099135}
R.~{Pries}.
\newblock Families of wildly ramified covers of curves.
\newblock {\em Amer. J. Math.}, 124(4):737--768, 2002.

\bibitem{MR2569747}
R.~{Pries}.
\newblock The {$p$}-torsion of curves with large {$p$}-rank.
\newblock {\em Int. J. Number Theory}, 5(6):1103--1116, 2009.

\bibitem{Seshadri1958-1959}
C.~S. Seshadri.
\newblock L'op{\'e}ration de {C}artier. {A}pplications.
\newblock {\em S{\'e}minaire Claude Chevalley}, 4:1--26, 1958-1959.

\bibitem{vanderGeer1999}
G.~van~der Geer.
\newblock {\em Cycles on the Moduli Space of Abelian Varieties. \upshape In:
  {Moduli of Curves and Abelian Varieties: The Dutch Intercity} {Seminar on
  Moduli}}, pages 65--89.
\newblock Vieweg+Teubner Verlag, Wiesbaden, 1999.

\bibitem{2018arXiv181204996Z}
Zijian {Zhou}.
\newblock {Ekedahl-Oort strata on the moduli space of curves of genus four}.
 arXiv:1812.04996, Dec 2018.

\end{thebibliography}
\bibliographystyle{plain}

\Addresses

\end{document}